\documentclass{article}
\usepackage{amssymb}
\usepackage[T2A]{fontenc} \usepackage[cp1251]{inputenc}
\usepackage[russian,ukrainian,english]{babel}
\def\finproof{\hbox{\vrule width1.0ex height1.ex}\vspace{2mm}}

\begin{document}

\centerline{\large On time transformations for differential
equations with state-dependent delay}

\medskip

\begin{center}
{\sc A.V. Rezounenko}

\smallskip

Kharkiv National University, Department of Mechanics and Mathematics,\\
4, Svobody Sqr., Kharkiv, 61022, Ukraine \\
{{\it E-mail:} rezounenko@yahoo.com}
\end{center}

\medskip

 Systems of differential equations with
state-dependent delay are considered. The delay dynamically depends
on the state i.e. is governed by an additional differential
equation. By applying the time transformations we arrive to constant
delay systems and compare the asymptotic properties of the original
and transformed systems.

\medskip


\centerline{\bf 1. Introduction}

\medskip

Delay differential equations successfully used to model and study a
number of applied problems in physics, biology, chemistry. Basic
mathematical theory of constant delay equations can be found in the
classical monographs \cite{Walther_book,Hale} and references
therein.

Differential equations with state-dependent delay (SDD) attracted
much attention during last decades and many deep results have been
obtained so far (see
\cite{Hartung-Krisztin-Walther-Wu-2006,Mallet-Paret,Walther_JDE-2003,Krisztin-2003}
and references therein). Such equations with discrete
state-dependent delays are always {\it nonlinear} by their nature.
As described in \cite{Hartung-Krisztin-Walther-Wu-2006}, this type
of delay brings additional difficulties in proving such basic
properties of solutions as uniqueness and continuous dependence on
initial data. The main approach to get the well-posed initial-value
problem is to restrict the set of initial functions and hence the
set of solutions to $C^1$ -functions
\cite{Hartung-Krisztin-Walther-Wu-2006}. For an alternative approach
where an additional condition on the SDD is used to get a well-posed
initial-value problem in the space of continuous functions see
\cite{Rezounenko_NA-2009,Rezounenko_JMAA-2012}. In this note we rely
on the classical approach \cite{Hartung-Krisztin-Walther-Wu-2006}
and compare the SDD problem with another one with a constant delay.
This constant delay problem is constructed by using the so-called
{\it time transformations}
\cite{Brunner-Maset-DCDS-2009,Brunner-Maset-CPAA-2010}. This
transformation could be applied to any particular solution along
which the deviating argument is monotone. To be assured 
that the monotonicity holds for all solutions we concentrate on the
system when the state-dependent delay is governed by an additional
differential equation and provide a sufficient condition for the
monotonicity of the deviating argument. This type of equations is
used to describe some models of population dynamics, see
\cite{Arino-Hadeler-Hbid-JDE-1998} and references therein. In
 \cite{Arino-Hadeler-Hbid-JDE-1998} one could also find motivations
 to study this type of SDD and  comparison with the frequently used case when SDD is presented by
 explicit or implicit functionals.

Our main goal in this note is to compare the asymptotic properties
of the SDD system with the corresponding ones of the system after
the time transformation. We try to find conditions which guarantee
that such properties as stability, boundedness and compactness of
the initial SDD problem survives under the time transformations.

\centerline{\bf 2. Time transformations}

\medskip

We study the following non-autonomous system with state-dependent
delay (see the autonomous case in
\cite{Arino-Hadeler-Hbid-JDE-1998})
\begin{eqnarray}
  \dot y(t) & = &f(t,y(t), y(t-\eta(t))), \quad t> t^0,\label{tt1}\\
  \dot \eta (t) & = &- \mu (\eta(t)-\widetilde \eta) + G(y(t)), \quad t> t^0,\label{tt2}
\end{eqnarray}
with initial data
\begin{eqnarray}
  y(t) & = & g(t), \quad t\in [t^0-h,t^0], \label{tt3}\\
  \eta (t^0) & =&  \eta^0.\label{tt4}
\end{eqnarray}

Here $y\in R^m, \eta \in R,  \eta^0,\widetilde{\eta}>0$. The
function $\eta$ is a state-dependent delay since it is a solution of
the equation (\ref{tt2}) where there is a dependence on $y$.

In the sequel we will denote $h\equiv 2\widetilde{\eta}$ and also
$X\equiv C^1([-h,0];R^m)\times R$ with the natural norm.

\medskip

{\bf Lemma~1.}  {\it Let $f$ and $G$ be Lipschitz and $|G(y)|\le \mu
\widetilde{\eta}$ for all $y\in R^m$. Then for any $g\in
C^1([t^0-h,t^0];R^m)$, $\eta(t^0)\equiv \eta^0\in [0,2\widetilde{\eta}] $
the system (\ref{tt1})-(\ref{tt4}) has a unique global solution such
that $\eta(t)\in [0,2\widetilde{\eta}] $ for all $t\ge t^0$. The
solution continuously depends on initial data $(g;\eta^0)$.
}%

\smallskip

{\it Proof of lemma~1}. The existence of solutions is simple since
the righthand sides of equations (\ref{tt1}), (\ref{tt2}) are
continuous. Solutions are global due to Lipschitz properties of $f$
and $G$. The uniqueness follows from the well-known results on the
state-dependent delay equations (see e.g.
~\cite{Hartung-Krisztin-Walther-Wu-2006}) since we consider
Lipschitz initial function $g$.

Using the property $|G(y)|\le \mu \widetilde{\eta}$, one can easily
show that (for any $y$) any solution of (\ref{tt2}), (\ref{tt4})
satisfies $\eta(t)\in [0,2\widetilde{\eta}] $ provided
$\eta(t^0)\equiv \eta^0\in
[0,2\widetilde{\eta}] $. 

Now we show the continuous dependence on initial data. For the
simplicity of presentation we put $t^0=0$. Let us consider a pair
$(\bar g;\bar\eta^0)\in C^1([-h,0];R^m)\times [0,h] \subset X$ and
an arbitrary sequence $(g^n;\eta^{0,n})$ such that $
||(g^n;\eta^{0,n}) - (\bar g;\bar\eta^0)||_X \to 0$ as $n\to
\infty$.

We rewrite the system  (\ref{tt1})-(\ref{tt4}) in the integral form

\begin{equation}\label{tt5}
y^n(t)= g^n(0)+\int^t_0 f(s,y^n(s), y^n(s-\eta^n(s))\, ds,
\end{equation}
\begin{equation}\label{tt6}
  \eta^n(t)-\widetilde \eta= e^{-\mu t}(\eta^{0,n}-\widetilde{\eta}) +\int^t_0 e^{-\mu (t-s)} G(y^n(s))\, ds.
\end{equation}
Similar equations are for initial data $(\bar g;\bar\eta^0)$. For
the differences of solutions, using the Lipschitz properties of $f$
and $G$ (the corresponding Lipschitz constants are denoted by $L_f$
and $L_G$) we have
$$ |y^n(t)-\bar y(t)| \le |g^n(0)-\bar g(0)| $$
$$+ L_f \int^t_0 \left\{ |y^n(s-\eta^n(s))-\bar y(s-\eta^n(s))| +
|\bar y (s-\eta^n(s))-\bar y (s-\bar\eta(s)) | + |y^n(s)-\bar y(s) |
\right\}\, ds,
$$
$$|\eta^n(t)-\bar\eta(t)| \le e^{-\mu t}|\eta^{0,n}- \bar\eta^0| +
L_G \int^t_0 |y^n(s)-\bar y(s) | \, ds.
$$
Let us fix any $T>0$. Since all solutions are $C^1$ in time (see
e.g. \cite{Hartung-Krisztin-Walther-Wu-2006})) for $t>0$ we can
denote by $L_{\bar y, T}$ the Lipschitz constant of the solution
$\bar y(t), t\in [-h,T]$. Hence $|\bar y (s-\eta^n(s))-\bar y
(s-\bar\eta(s)) | \le L_{\bar y, T} |\eta^n(s)-\bar\eta(s)|$ for all
$s\in [0,t]\subset [0,T]$. Denoting for short $\beta^n(t)\equiv
\max_{\tau\in [0,t]} \{ |y^n(\tau)-\bar y(\tau)| +
|\eta^n(\tau)-\bar\eta(\tau)| \}$ and $C^T\equiv \max \{ 2L_f
+L_G;L_{\bar y, T} \}$ we obtain
$$0\le \beta^n(t) \le \beta^n(0) + L_f \max_{\tau\in [-h,0]}|g^n(\tau)-\bar
g(\tau)| + C^T\int^t_0 \beta^n(s)\, ds.
$$
We apply the Gronwall's inequality and get for all $t\in [0,T]$
$$\max_{\tau\in [0,t]} \{ |y^n(\tau)-\bar y(\tau)| + |\eta^n(\tau)-\bar\eta(\tau)|
\} $$ $$ \le \left\{ (1+L_f)\max_{\tau\in [-h,0]}|g^n(\tau)-\bar
g(\tau)| + |\eta^{0,n}- \bar\eta^0| \right\} \exp( t\cdot \max \{
2L_f +L_G;L_{\bar y, T} \} ).
$$
The last convergence and equations (\ref{tt1}), (\ref{tt2}) give a
similar estimate (convergence) for the time derivatives, so we
finally get
$$ |y^n-\bar y|_{C^1[0,t]} + |\eta^n-\bar\eta|_{C^1[0,t]}
 \to 0 \quad \hbox{ as } \quad  n\to + \infty.
$$

This gives the continuous dependence on initial data and completes
the proof of lemma~1. \finproof 

\medskip

For any solution $(y;\eta)$ of the system (\ref{tt1})-(\ref{tt4}) we
call the function $\sigma$ given by
\begin{equation}\label{tt7}
    \sigma(t)= t - \eta(t), \quad t\ge t^0
\end{equation}
the {\sl deviating argument} of  $(y;\eta)$.

Our goal is to investigate properties connected to the time
transformation approach introduced in
\cite{Brunner-Maset-DCDS-2009,Brunner-Maset-CPAA-2010}. We are going
to use a function $t=\alpha(s)$ called {\sl time
transformation}~\cite{Brunner-Maset-CPAA-2010} to convert a
particular solution $(y;\eta)$ of the system (\ref{tt1})-(\ref{tt4})
into a solution $(z;\chi;\alpha)$ of {\sl constant} delay system
\begin{equation}\label{tt8}
   \left\{ \begin{array}{ccl}
           \dot z(s)   & = & f(\alpha(s),z(s),z(s-h))\, \dot \alpha(s), \quad s\ge s^0, \\
            z(s) & = & g(\omega(s)), \quad s\in [s^0-h,s^0], \\
            \dot \chi(s) &= & - \mu (\chi(s)-\widetilde \eta )\, \dot \alpha(s) +
            G(z(s))\, \dot \alpha(s),\\
            \chi(s^0) & = & \eta^0,
           \end{array}
   \right.
\end{equation}
where $\alpha$ satisfies the algebraic equation
\begin{equation}\label{tt9}
   \left\{ \begin{array}{ccl}
           \alpha(s)- \chi(s) & =   & \alpha(s-h), \quad s\ge s^0, \\
            \alpha(s) & = & \omega(s), \quad s\in [s^0-h,s^0], \\
           \end{array}
   \right.
\end{equation}
Here $\omega : [s^0-h,s^0]\to R$ is an arbitrary $C^1$-function with
positive derivative 
and such that $\omega(s^0-h)=\omega(s^0) - \eta^0 <t^0,
\omega(s^0)=t^0.$

\smallskip

{\bf Remark.} {\it We notice that equation (\ref{tt9}) is different from
the corresponding rules used in \cite{Brunner-Maset-CPAA-2010} and
\cite{Brunner-Maset-DCDS-2009} since
here we  have no a given {\tt lag function}.
}%

The time transformation $\alpha$ is constructed step by step (see
(\ref{tt9}), (\ref{tt7})) by the rule
\begin{equation}\label{tt10}
\alpha(s)= \sigma^{-1}(\alpha(s-h)), \quad s\in [s^0+(k-1)h,s^0+kh],
\quad k=0,1,2... .
\end{equation}
Here we used (see (\ref{tt7})) $ \sigma(\alpha(s))=\alpha(s) -
\eta(\alpha(s))=\alpha(s) - \chi(s) = \alpha(s-h)$, since
$\chi(s) = \eta(\alpha(s))$.

It is clear that one needs the invertibility of $\sigma$ to define
$\alpha$. In \cite{Brunner-Maset-CPAA-2010} the rule (\ref{tt10})
was used {\sl assuming} that $\dot \sigma (t) >0$ (or $\dot \sigma
(t) <0$). More precisely, it was used for those solutions along
which $\dot \sigma (t) >0$ (or $\dot \sigma (t) <0$). In our study
we can give a simple condition which {\sl guarantees} that along
{\sl all} solutions we have $\dot \sigma (t) >0$ and hence $\sigma$
is always invertible. Such a simple condition is
$2\mu\widetilde{\eta}<1$. It is easy to see using (\ref{tt2}) that
in this case $|\dot \eta(t)| \le \mu |\eta-\widetilde{\eta}| +
|G(y(t))| \le \mu\widetilde{\eta} + \mu\widetilde{\eta} <1$. Here we
used the assumption $|G(y)|\le \mu \widetilde{\eta}$. Now
(\ref{tt7}) implies $\dot \sigma (t) = 1 - \dot \eta(t) >0$.

\smallskip

{\bf Remark.} {\it Taking into account that
the state-dependent delay $\eta$ takes values in
$[0,2\widetilde{\eta}]=[0,h]$, one could say that the assumption
$2\mu\widetilde{\eta}<1$ means the "slow changing of delay" in the
range $[0,2\widetilde{\eta}]$.}

\smallskip

It is important that if $\dot \sigma (t) >0$, then $\dot \alpha (s)
>0$ since on the initial time segment $\dot \alpha (s)=\dot \omega
(s)>0, s\in [s^0-h,s^0]$ and $\dot \alpha (s)=\frac{d
\sigma^{-1}(\alpha(s-h))}{ds}=\frac{1}{\dot \sigma(\alpha(s-h))}
\dot \alpha(s-h)=\frac{\dot \alpha(s-h)}{\dot
\sigma(\sigma(\alpha(s)))}>0$ step by step (see also
\cite[p.28]{Brunner-Maset-CPAA-2010}). Here we used
$\sigma(\alpha(s))=\alpha(s-h)$ (see (\ref{tt9}), (\ref{tt10})).

\smallskip

By construction (see \cite[propositions
1,2]{Brunner-Maset-CPAA-2010}) the connection between a solution
$(y,\eta)$ of (\ref{tt1})-(\ref{tt4}) and the corresponding solution
$(z,\chi,\alpha)$ of (\ref{tt8})-(\ref{tt9}) is given by
\begin{equation}\label{tt11}
y(t)= z(\alpha^{-1}(t)), \, t\in [t^0-h,+\infty), \, z(s) =
y(\alpha(s)),  \, \chi(s)=\eta(\alpha(s)),\, s\in [s^0-h,+\infty),
\, t=\alpha(s),\, t^0=\alpha(s^0).
\end{equation}

\medskip

\medskip

{\bf Lemma~2.}  {\it Let $f$ and $G$ be as in lemma~1. Consider a
sequence $\{(g^n;\eta^{0,n})\}$ such that $ ||(g^n;\eta^{0,n}) -
(\bar g;\bar\eta^0)||_{C^1([t^0-h,t^0];R^m)\times R} \to 0$ and a
sequence $\{\omega^n\}$ such that $||\omega^n-
\bar\omega||_{C^1[t^0-h,t^0]}\to 0 $ as $n\to \infty$. Then for any
$S>0$ the sequence of time transformations $\alpha^n$ uniformly
converges to $\bar\alpha$ i.e. $\max_{s\in [s^0,s^0+S]} |\alpha^n(s)
- \bar\alpha (s)| \to 0$ as $n\to \infty$
}%

\smallskip

{\it Proof of lemma~2}. First, using  (\ref{tt7}), we have from lemma~1 that
$\max_{t\in [t^0,t^0+T]}|\sigma^n(t)-\bar\sigma(t)|\to 0$ as $n\to \infty$.
This convergence and (\ref{tt2}) imply that $||\sigma^n-\bar\sigma||_{C^1[t^0,t^0+T]}\to 0$.

Since ${d\over dt}\bar\sigma(t)>0$ on $[t^0,t^0+T]$, then there is a $\bar\delta>0$
such that ${d\over dt}\bar\sigma(t)\ge 2\bar\delta>0$ for all $t\in [t^0,t^0+T]$.
This and $||\sigma^n-\bar\sigma||_{C^1[t^0,t^0+T]}\to 0$ give
\begin{equation}\label{tt20}
{d\over dt} \sigma^n(t)\ge \bar\delta>0 \hbox{ for all } t\in [t^0,t^0+T]\hbox{  and } n>n_1.
\end{equation}

Using the definition of $\alpha$ (see (\ref{tt10})) and the
convergence $\omega^n\to \bar\omega $ as $n\to \infty$ we only need
to show the property $\max_{\tau\in [0,T]} |(\sigma^n)^{-1}(\tau) -
(\bar\sigma)^{-1}(\tau)| \to 0$ as $n\to \infty$. Let us denote by
$\gamma^n(s) =(\sigma^n)^{-1}(s), \bar \gamma
(s)=(\bar\sigma)^{-1}(s)$. Assume opposite i.e. $\gamma^n(s)
\not\rightrightarrows  \bar \gamma (s)$ on some $[s^0,s^0+S]$. Hence
$$ \exists \varepsilon_0>0\, : \, \forall
N\in {\bf N}\quad \exists\, n_N\ge N \quad\exists\, s_{n_N}\in
[s^0,s^0+S] \, :\, |\gamma^n(s_{n_N}) - \bar \gamma (s_{n_N}) | \ge
\varepsilon_0.
$$
Considering $N=1,2,...$ we get two sequences $\{ n_k\}^\infty_{k=1}$
and $\{ s_{n_k}\}^\infty_{k=1}\subset [s^0,s^0+S]$ such that
\begin{equation}\label{tt21}
|\gamma^{n_k}(s_{n_k}) - \bar \gamma (s_{n_k}) | \ge \varepsilon_0.
\end{equation}
Since $[s^0,s^0+S]$ is compact we have $\hat s\in [s^0,s^0+S]$ and a
subsequence again denoted by $\{ n_k\}^\infty_{k=1}$ such that
$s_{n_k} \to \hat s\in [s^0,s^0+S]$. We can write
$$ \gamma^{n_k}(s_{n_k}) - \bar \gamma (s_{n_k}) =
\left( \gamma^{n_k}(s_{n_k}) - \gamma^{n_k}(\hat s)\right) +
\left(\gamma^{n_k}(\hat s) - \bar \gamma (\hat s)\right) +
\left(\bar \gamma (\hat s) - \bar \gamma (s_{n_k})\right)
$$
The last term vanishes due the continuity of $\bar \gamma$, the second one due to
the point-wise convergence $\gamma^n(s) \to \bar \gamma (s)$ for all
$s\in [s^0,s^0+S]$.
Hence the only possibility to satisfy (\ref{tt21}) is that there is an integer
$k_1$ such that for all $k\ge k_1$ one has
$$|\gamma^{n_k}(s_{n_k}) - \gamma^{n_k}(\hat s)| >\varepsilon_0/2.$$
The last property together with $s_{n_k} \to \hat s$ and differentiability of all
 $\gamma^{n_k}$ imply that the derivatives ${d\over ds}\gamma^{n_k}$ are unbounded
in a neighborhood of $\hat s$. This contradicts the property
${d\over dt} \sigma^n(t) \ge \bar\delta>0$ (see (\ref{tt20})) since
${d\over ds} \gamma^n(s)= {1\over {d\over dt}\sigma^n(t) } $.
 The proof of lemma~2 is complete. \finproof

\medskip

Now combining lemmata 1 and 2 (and property $\mu \widetilde{\eta}<
{1\over 2}$) we can formulate the first result on the continuous
dependence of the time transformation on initial data.

\medskip

{\bf Theorem 1.} {\it Let $f$ and $G$ be Lipschitz and $|G(y)|\le
\mu \widetilde{\eta}< {1\over 2}$ for all $y\in R^m$. Consider a
sequence $\{(g^n;\eta^{0,n})\}$ such that $ ||(g^n;\eta^{0,n}) -
(\bar g;\bar\eta^0)||_{C^1([t^0-h,t^0];R^m)\times R} \to 0$ and a
sequence $\{\omega^n\}$ such that $||\omega^n-
\bar\omega||_{C^1[t^0-h,t^0]}\to 0 $ as $n\to \infty$. Then the time
transformation gives the sequence of the corresponding solutions $\{
(z^n,\chi^n,\alpha^n)\}$ of the constant delay system
(\ref{tt8})-(\ref{tt9}) (see (\ref{tt11})) such that  for any $S>0$
one has
$$ \max_{s\in [s^0,s^0+S]} \left\{ ||z^n(s)-\bar z(s)||  +  |\chi^n(s)
- \bar\chi (s)|+  |\alpha^n(s) - \bar\alpha (s)|\right\} \to 0
\,\,\hbox{ as }\,\, n\to \infty.
$$
and $\dot \alpha^n(s)>0$ for $s\in [s^0,s^0+S]$.}

\smallskip

{\bf Remark.} {\it  We should notice that (\ref{tt8}), (\ref{tt9})
is the system of coupled differential and algebraic equations. It is
necessary to comment on how to solve it. The way is different from
the one of \cite{Brunner-Maset-CPAA-2010,Brunner-Maset-DCDS-2009}
since we have no a given lag function (c.f. \cite[Section
2.1]{Brunner-Maset-CPAA-2010}). Using (\ref{tt9}), we write for
$s\in [0,h]$ (and than continue step by step with the step $h$):
$\alpha(s)=\chi (s) + \omega(s-h)$. Then we substitute it into the
differential equation for $\chi$ in (\ref{tt8}) to get $\dot \chi(s)
= - \mu (\chi(s)-\widetilde \eta )\, (\dot \chi (s) + \dot
\omega(s-h)) + G(z(s))\, (\dot \chi (s) + \dot \omega(s-h)).$ Hence
$\dot \chi(s) \left[1 + \mu (\chi(s)-\widetilde \eta ) -
G(z(s))\right] = \left\{- \mu (\chi(s)-\widetilde \eta )\, ) +
G(z(s)\right\}\,  \dot \omega(s-h).$ We remind that the assumption
$2\mu \widetilde \eta <1$ implies $|\mu (\chi(s)-\widetilde \eta ) -
G(z(s)) |<1$.  It gives
\begin{equation}\label{tt25}
 \dot \chi(s)  = \left\{- \mu (\chi(s)-\widetilde \eta )\, ) +
G(z(s)\right\}\,  \dot \omega(s-h) \left[1 + \mu (\chi(s)-\widetilde
\eta ) - G(z(s))\right]^{-1}.
\end{equation}
Now to rewrite the first equation in (\ref{tt8}) we use again
$\alpha(s)=\chi (s) + \omega(s-h)$ (and $\dot \alpha(s)=\dot \chi
(s) + \dot \omega(s-h)$), substitute it in (\ref{tt8}) and use
(\ref{tt25}). It gives $\dot z(s)   =f(\alpha(s),z(s),z(s-h))\, \dot
\alpha(s) = f(\chi (s) + \omega(s-h),z(s),z(s-h))\, (\dot \chi (s) +
\dot \omega(s-h)) = f(\chi (s) + \omega(s-h),z(s),z(s-h))\, \left(
\frac{- \mu (\chi(s)-\widetilde \eta )\, ) + G(z(s)}{1+  \mu
(\chi(s)-\widetilde \eta )\, ) - G(z(s)} +1 \right) \cdot \dot
\omega(s-h)=  f(\chi (s) + \omega(s-h),z(s),z(s-h))\, \left[1 + \mu
(\chi(s)-\widetilde \eta ) - G(z(s))\right]^{-1}\cdot \dot
\omega(s-h).$

Hence (\ref{tt8}), (\ref{tt9}) can be rewritten for $s\in [0,h]$
(the first step) as
\begin{equation}\label{tt26}
   \left\{ \begin{array}{ccl}
           \dot z(s)   & = & f(\chi (s) + \omega(s-h),z(s),z(s-h))\,
           \left[1 + \mu (\chi(s)-\widetilde \eta ) - G(z(s))\right]^{-1}
           \cdot \dot \omega(s-h), \quad s\in [0,h],\\ 
            z(s) & = & g(\omega(s)), \quad s\in [0,h],\\ 
            \dot \chi(s) &= & \left\{- \mu (\chi(s)-\widetilde \eta )\, ) +
G(z(s)\right\}\,  \dot \omega(s-h) \left[1 + \mu (\chi(s)-\widetilde
\eta ) - G(z(s))\right]^{-1},\\
            \chi(0) & = & \eta^0,
           \end{array}
   \right.
\end{equation}
It easy to see that the last system gives solution $(z(s),\chi(s))$
for $s\in [0,h]$.  Then the time transformation $\alpha$ is found by
$\alpha(s)=\chi (s) + \omega(s-h)$. In general it reads as
\begin{equation}\label{tt27}
   \left\{ \begin{array}{ccl}
           \dot z(s)   & = & f(\chi (s) + \alpha(s-h),z(s),z(s-h))\,
           \left[1 + \mu (\chi(s)-\widetilde \eta ) - G(z(s))\right]^{-1}
\cdot \dot \alpha(s-h), \quad s\ge s^0, \\
            z(s) & = & g(\omega(s)), \quad s\in [s^0-h,s^0], \\
            \dot \chi(s) &= & \left\{- \mu (\chi(s)-\widetilde \eta )\, ) +
G(z(s)\right\}\,  \left[1 + \mu (\chi(s)-\widetilde
\eta ) - G(z(s))\right]^{-1} \dot \alpha(s-h) ,\quad s\ge s^0,\\
            \chi(s^0) & = & \eta^0,
           \end{array}
   \right.
\end{equation}
and can be solved step by step.
}
\bigskip

\centerline{\bf 3. Connection between asymptotic properties of
systems (\ref{tt1})-(\ref{tt4}) and (\ref{tt8})-(\ref{tt9}).}

\bigskip

In this section we discuss how to determine if some qualitative
properties of solutions 
of the initial state-dependent delay system
 (\ref{tt1})-(\ref{tt4}) survive the time transformation i.e. still be valid for the corresponding
 solutions 
 of constant delay system (\ref{tt8})-(\ref{tt9}). We are also
 interested to connect the known properties of solutions of
 (\ref{tt8})-(\ref{tt9}) with the ones of  (\ref{tt1})-(\ref{tt4}).

\smallskip

Let us start with the discussion of the property of the exponential
stability. For the simplicity of presentation we assume that the
zero function $y(t)\equiv 0$ is a solution of
(\ref{tt1})-(\ref{tt4}). Hence, by (\ref{tt11}), $z(s)\equiv 0$ will
be also a solution of (\ref{tt8})-(\ref{tt9}).

Assume one has
$$||z(s)|| \le D^0 e^{-D^1 (s-s^0)}
||z_{s^0}||_{C([-h,0];R^m)}, s\ge s^0, D^0,D^1>0.
$$
 Hence, by
(\ref{tt11}), we get $$||y(t)||\le D^0 e^{-D^1 (\alpha^{-1}(t)-s^0)}
||z_{s^0}||_{C([-h,0];R^m)}= D^0 e^{-D^2 (t-t^0)}
||z_{s^0}||_{C([-h,0];R^m)} e^{D^2 (t-t^0) -D^1
(\alpha^{-1}(t)-s^0)}.$$ It is easy to see that if (and only if)
$e^{D^2  t -D^1 \alpha^{-1}(t)}$ is bounded, then we have
 $$||y(t)||\le D^3 e^{-D^2 (t-t^0)}
||z_{s^0}||_{C([-h,0];R^m)}.
$$
The above considerations show that the exponential stability of the
zero solution $z(s)\equiv 0$ implies the exponential stability of
the zero solution $y(t)\equiv 0$ provided there are positive
constants $D^1,D^2$ such that ${D^2 t -D^1 \alpha^{-1}(t)}\le 0 $
for all $t\ge t^0$. Similar estimates give the inverse implication
i.e. the exponential stability of the zero solution $y(t)\equiv 0$
implies the exponential stability of the zero solution $z(s)\equiv
0$ provided there are positive constants $C^1,C^2$ such that ${C^2 s
-C^1 \alpha (s)}\le 0 $ for all $s\ge s^0$. Since by (\ref{tt11})
$t=\alpha(s)$, we arrive to the following 

\smallskip

{\bf Definition 1.} {\it We say that "$s$-time" and "$t$-time" are
{\sl equivalent} if there are 
constants $A^1>0, A^2>0, B^1\in R, B^2\in R$ such that $A^1 t
+B^1\le s \le A^2 t+B^2$. }

\smallskip

{\bf Remark A.} {\it It is evident that in this case we also have $
{1\over A^2} s -{B^2\over A^2}\le t \le {1\over A^1} s -{B^1\over
A^1}$.}

\smallskip

We saw that in the case of the equivalent "$s$-time" and "$t$-time"
the property of the exponential stability survives under the time
transformation. Another consequence of the equivalence is that $t\to
+\infty$ if and only if $s\to +\infty$, which is clearly important
for the study of long-time asymptotic behavior of solutions.

\smallskip

The last result suggests to study the notion of time-equivalence in
detail. Let us try to find if in our case we have the equivalence.
The rule (\ref{tt10}) and (\ref{tt7}) show that we need to analyse
the function $\sigma^{-1}$. Using the property $\sigma(t) \ge t-h$
(bounded delay) and invertibility of $\sigma$, we get
$\sigma^{-1}(\tau) \le \tau +h$. Hence, by (\ref{tt10}), one has
$\alpha(s)= \sigma^{-1}(\alpha(s-h)) \le \alpha(s-h) + h.$
Particularly, $\alpha(h)\le \alpha(0) +h, \alpha(2h)\le \alpha(h) +
h \le \alpha(0) +2h, etc.$ Hence the property $\alpha(kh)\le
\alpha(0) +kh$ and the strict increase of $\alpha$ give the
following estimate
\begin{equation}\label{tt22}
\alpha (s) \le (\alpha(0) + h) + s.
\end{equation}
Since $\alpha(s)=t$, the estimate (\ref{tt22}) means the lower bound
in definition 1 and the upper bound in remark~A (with $A^1=1,
B^1=-(\alpha(0) + h) $).

The supplemented bounds in definition 1 and remark~A are less
obvious. For the moment we do not  claim that it is true in general
case, but present an additional assumption which guarantees the
bounds.

Let us assume that the value of delay is bounded from below by a
positive constant, say  $h_1>0 $. More precisely, $\eta(t)\ge h_1\in
(0, \widetilde{\eta}]$ for all $t\ge 0$. A sufficient condition for
the last property is $|G(y)|\le \mu |\widetilde{\eta}-h_1|$ for all
$y\in R^m$. Under the above condition we have $\sigma^{-1}(\tau) \ge
\tau +h_1$. By (\ref{tt10}), one has $\alpha(s)=
\sigma^{-1}(\alpha(s-h)) \ge \alpha(s-h) + h_1.$

Particularly, $\alpha(h)\ge \alpha(0) +h_1, \alpha(2h)\le \alpha(h) +
h_1 \le \alpha(0) +2h_1, etc.$ Hence the property $\alpha(kh)\le
\alpha(0) +kh_1$ and the strict increase of $\alpha$ give the
following estimate
\begin{equation}\label{tt23}
\alpha (s) \ge (\alpha(0) - h_1) + {h_1\over h}s.
\end{equation}
Combining (\ref{tt22}) and (\ref{tt23}) we conclude that "$s$-time"
and "$t$-time" are equivalent  (with $A^1=1, B^1=-(\alpha(0) + h),
A^2={h \over h_1},  B^2=- {h \over h_1} (\alpha(0) - h_1)$ in
definition 1 and remark~A).

Having the equivalence proved we can use it to compare the
asymptotic behavior of the corresponding dynamical systems
(process), constructed by solutions of systems before and after the
time transformations.


\smallskip

Let us consider an autonomous case of (\ref{tt1}) i.e.  the system
(c.f. (\ref{tt1})-(\ref{tt4}))
\begin{eqnarray}
  \dot y(t) & = &f^a(y(t), y(t-\eta(t))), \quad t> 0,\label{tt1auto}\\
  \dot \eta (t) & = &- \mu (\eta(t)-\widetilde \eta) + G(y(t)), \quad t> 0,\label{tt2auto}
\end{eqnarray}
with initial data
\begin{eqnarray}
  y(t) & = & g(t), \quad t\in [-h,0], \label{tt3auto}\\
  \eta (0) & =&  \eta^0.\label{tt4auto}
\end{eqnarray}

 We can restrict our study (using lemma 1) to the set of initial data
\begin{equation}\label{tt24}
X_{f^a}=\left\{ (g,\eta^0) \, | \,  \dot g(0) =
f^a(g(0),g(-\eta^0))\right\} \subset C^1([-h,0];R^m)\times [0,h].
\end{equation}
The set $X_{f^a}$ is an analogous to the {\it solution manifold}
used in \cite{Walther_JDE-2003} (see also
 \cite{Hartung-Krisztin-Walther-Wu-2006}). We notice that the reason
 for this restriction is to have $C^1$ smoothness of solution at
 zero i.e. $\dot y(0-)=\dot g(0-)= \dot y(0+)=f^a(g(0),g(-\eta^0))$. It is  easy to see that
 $X_{f^a}$ is invariant.

 We define the evolution operator $S^a (t) : X_{f^a} \to X_{f^a}$,  associated to the
 system (\ref{tt1auto})-(\ref{tt4auto}), by the formula $S^a(t) (g;\eta^0)=
(y_t;\eta(t))$, where $(y;\eta)$ is the unique solution of
(\ref{tt1auto})-(\ref{tt4auto}). It is easy to see that under our
assumptions the pair $(S^a, X_{f^a})$ constitutes a dynamical system
(in other words, the IVP (\ref{tt1auto})-(\ref{tt4auto}) is
well-posed in $X_{f^a}$). For more definitions and details on
dynamical systems see e.g.
\cite{Hale_AMS-1988_book,Temam_Springer-1988_book,Chueshov_Acta-1999_book}.

\smallskip

To discuss the 
 properties of solutions to the {\it
non-autonomous} system  (\ref{tt8}), (\ref{tt9}), let us remind the
following definition from \cite[pp. 112-119]{Vishik-book-1992}. Let
$E$ be a Banach space. Consider the two-parameter family of maps $\{
U(t,\tau) \}$, $U(t,\tau) :E \to E$, parameters $\tau\in R, t\ge
\tau$.

\smallskip

{\bf Definition 2.} \cite[p.113]{Vishik-book-1992}. {\it A family of
maps $\{ U(t,\tau) \}$ is called a {\tt process} on $E$ if
\par (i) $U(\tau,\tau)=I$=indentity;
\par (ii) $U(t,s) \circ U(s,\tau)=U(t,\tau)$ for all $t\ge s\ge \tau
\in R.$
}%

\smallskip

Since (\ref{tt8}), (\ref{tt9}) has constant delay $h$ we have no
need to restrict our study to a solution manifold. We define
$$E=
C([-h,0];R^m)\times [0,h] \times $$ $$\times \{ \omega \in
C^1[-h,0]\, | \, \dot \omega (\cdot) >0, \, \omega(-h)=\omega(0) -
\eta^0;\, \dot \omega(0) \left[1 + \mu (\eta^0-\widetilde \eta ) -
G(\varphi(0))\right] =\dot \omega(-h)\}$$
 and define $U(s,s^0) (\varphi;\eta^0;\omega)=
(z_s;\chi(s);\alpha_s),$ where $(z;\chi;\alpha)$ is the unique
solution of (\ref{tt8}), (\ref{tt9}) with initial data
$(\varphi;\eta^0;\omega),$ here $ z_{s^0}=\varphi $.

We notice that to come back to the  system
 (\ref{tt1auto})-(\ref{tt4auto}), using a particular solution of (\ref{tt8}),
 (\ref{tt9}), one restores function $g$ as follows $g(t) =
 z_{s^0}(\alpha^{-1}(t))$ for $t\in [t^0-h,t^0]$.




\smallskip

Let us continue to discuss which asymptotic properties of
(\ref{tt1})-(\ref{tt4}) survive the time transformation i.e. still
be valid for the corresponding
 solutions 
 of constant delay system (\ref{tt8})-(\ref{tt9}).

One of the important properties of dynamical systems and processes
are boundedness of solutions and their compactness (asymptotic
compactness) see e.g.
\cite{Hale_AMS-1988_book,Temam_Springer-1988_book,Chueshov_Acta-1999_book}.
Below we always assume that $t\to +\infty$ if and only if $s\to
+\infty$, which is true, for example, in the case of the equivalence
of $t-time$ and $s-time$. One can easily see, by (\ref{tt11}), that
$||y(t)|| \le C, t\ge t_1$ is equivalent to  $||z(s)|| \le C, s\ge
s_1$. Hence, the existence of a bounded absorbing set for $z_s$
-coordinate is equivalent to the existence of a bounded absorbing
set for $y_t$ -coordinate, both in the space  $C([-h,0];R^m)$. To go
further, let us discuss the following additional assumptions on the
time transformation $\alpha$:
\par $ {\bf (A1)}\qquad \exists C^{1,\alpha}>0, \forall s\ge s_1
\Rightarrow \dot \alpha(s) \le C^{1,\alpha}.$
\par $ {\bf (A1^\prime)}\qquad \exists C^{2,\alpha}>0, \forall t\ge
t_1
\Rightarrow \dot \alpha^{-1}(t) \le C^{2,\alpha}.$
\par $ {\bf (A2)}\qquad \alpha \hbox{ is uniformly continuous on }
[s_1,+\infty). $
\par $ {\bf (A2^\prime)}\qquad \alpha^{-1} \hbox{ is uniformly continuous on }
[t_1,+\infty). $
\par $ {\bf (A3)}\qquad \dot \alpha \hbox{ is uniformly continuous on }
[s_1,+\infty). $
\par $ {\bf (A3^\prime)}\qquad \dot \alpha^{-1} \hbox{ is uniformly continuous on }
[t_1,+\infty). $
\smallskip

Using (\ref{tt11}), we  have $\dot z(s)=\dot y(\alpha(s))\, \dot
\alpha(s)$. Hence, the existence of a bounded absorbing set for the
$y_t$ -coordinate in the space $C^1([-h,0];R^m)$ implies the
existence of a bounded absorbing set for the $z_s$ -coordinate,
provided (A1) is satisfied. The inverse implication is valid
provided $(A1^\prime)$ is satisfied. For the system
(\ref{tt1})-(\ref{tt4}) the existence of a bounded absorbing set
means it is dissipative (for more details on this property see e.g.
\cite{Hale_AMS-1988_book,Temam_Springer-1988_book,Chueshov_Acta-1999_book}).
Let us now assume that  for $t\ge t_1$ the $y_t$ -coordinate belongs
to a (pre-) compact set in the space $C([-h,0];R^m)$. By the
Arzela–Ascoli theorem theorem, the family $\{ y_t \}_{t\ge t_1}$ is
uniformly bounded and equicontinuous. Using (\ref{tt11}), we  have
$|z(s^1)-z(s^2)| = |y(\alpha(s^1)) - y(\alpha(s^2))|$. This and the
above discussion show that for $s\ge s_1$ the $z_s$ -coordinate
belongs to a (pre-) compact set in the space $C([-h,0];R^m)$,
provided (A2) is satisfied. The inverse implication is valid
provided $(A2^\prime)$ is satisfied. The similar considerations in
$C^1([-h,0];R^m)$ need the estimate $|\dot z(s^1)-\dot z(s^2)|=|\dot
y(\alpha(s^1))\, \dot \alpha(s^1) - \dot y(\alpha(s^2))\, \dot
\alpha(s^2)| \le |\dot y(\alpha(s^1)) - \dot y(\alpha(s^2)) |\,
|\dot \alpha(s^1)|\, + |\dot y(\alpha(s^2))|\, |\dot \alpha(s^1) -
\dot \alpha(s^2) |.$ We see that if for $t\ge t_1$ the $y_t$
-coordinate belongs to a (pre-) compact set in the space
$C^1([-h,0];R^m)$, then for $s\ge s_1$ the $z_s$ -coordinate belongs
to a (pre-) compact set in the space $C^1([-h,0];R^m)$, provided
(A1)-(A3) are satisfied. The inverse implication is valid provided
$(A1^\prime)$ - $(A3^\prime)$ are satisfied. Particularly, we shown
that assumptions (A1)-(A3) and $(A1^\prime)$ - $(A3^\prime)$ connect
asymptotic properties of the dynamical system $(S^a(t), X_{f^a})$
and the process $U(s,\tau) :E \to E$.

\smallskip

{\bf Remarks.} {\it 1. Discussing assumptions (A1)-(A3), one could
think that the family $\{\alpha_s\}_{s\ge s_1}$ may belong to a
(pre-) compact set in the space $C[-h,0]$ or even $C^1[-h,0]$, but
it is never true since $\alpha(s)=t$ is time, which is naturally
unbounded.
\par 2. One can see that $(A1)$ implies $(A2)$, but $(A2)
\not\Rightarrow (A1)$, similarly $(A1^\prime) \Rightarrow
(A2^\prime),$ but $(A2^\prime) \not\Rightarrow (A1^\prime)$.
\par 3. We notice that $(A1)$ gives $\alpha (s) \le C^{1,\alpha} s +k^1$ and
similarly $(A1^\prime) \Rightarrow \alpha^{-1} (t) \le C^{2,\alpha}
t +k^2$. Using definition~1, we see that $(A1), (A1^\prime)$ imply
the equivalence of "s-time" and "t-time".
}%

\medskip

{\bf Acknowledgments}. The author thanks S. Maset for useful
comments on an earlier version of the manuscript.

%
%



\end{document}